\documentclass{amsart}

\usepackage{amssymb,amsmath, amsfonts}

\newtheorem{theorem}{Theorem}[section]

\newtheorem{corollary}[theorem]{Corollary}
\newtheorem{proposition}[theorem]{Proposition}

\theoremstyle{definition}

\theoremstyle{remark}
\newtheorem{remark}[theorem]{Remark}

\numberwithin{equation}{section}

\newcount\quantno
\everydisplay{\quantno=0}\everycr{\quantno=0}
\newcommand\quant{\advance\quantno by1
                      \ifnum\quantno=1\qquad\else\quad\fi\forall }
\newcommand\itemno[1]{(\romannumeral #1)}

\newcommand\rest[1]{\kern-.1em
          \lower.5ex\hbox{$\scriptstyle #1$}\kern.05em}

\renewcommand\mod[1]{\left\vert{#1}\right\vert}
\newcommand\bigmod[1]{\bigl\vert{#1}\bigr|}
\newcommand\Bigmod[1]{\Bigl\vert{#1}\Bigr|}

\newcommand\norm[2]{{\Vert{#1}\Vert_{#2}}}

\newcommand\prodo[2]{\left\langle#1,#2\right\rangle}

\newcommand\wrt{\,\text{\rm d}}

\newcommand\BC{\mathbb{C}}

\newcommand\BN{\mathbb{N}}
\newcommand\BR{\mathbb{R}}

\newcommand\BZ{\mathbb{Z}}

\newcommand\cA{\mathcal{A}}

\newcommand\cM{\mathcal{M}}

\newcommand\cS{\mathcal{S}}

\newcommand\be{\beta}
 
\newcommand\de{\delta}
 \newcommand\vep{\varepsilon}

\newcommand\la{\lambda} 
 \newcommand\Om{\Omega} 
\newcommand\si{\sigma}

\newcommand\vp{\varphi}

\newcommand\funnyk{k\hbox to 0pt{\hss\phantom{g}}}

\newcommand\lu[1]{L^1(#1)}

\newcommand\laq[1]{L^q(#1)}
\newcommand\laqO[1]{L_0^q(#1)}
\newcommand\laqcO[1]{L_{c,0}^q(#1)}
\newcommand\ld[1]{L^2(#1)}
\newcommand\ldO[1]{L^2_0(#1)}
\newcommand\ldcO[1]{L^2_{c,0}(#1)}
\newcommand\ldloc[1]{L^2_{\mathrm{loc}}(#1)}

\newcommand\ly[1]{L^\infty(#1)}

\newcommand\hu[1]{H^1(#1)}

\newcommand\hudfin[1]{H_{\mathrm{fin}}^{1,2}(#1)}
\newcommand\huqfin[1]{H_{\mathrm{fin}}^{1,q}(#1)}
\newcommand\huinftyfin[1]{H_{\mathrm{fin}}^{1,\infty}(#1)}

\newcommand\wt{\widetilde}
\newcommand\whH{\widehat{\phantom{G}}\hbox to 0pt{\hss $H$}}

\newcommand\emspace{\hbox to 6pt{\hss}}

\newcommand\rmi{\hbox{\rm (i)}}
\newcommand\rmii{\hbox{\rm (ii)}}

\newcommand\diam[1]{\mathrm{diam}#1}

\begin{document}

\title{On the $H^1$--$L^1$ boundedness of operators}

\author{Stefano Meda}
\address{Dipartimento di Matematica e Applicazioni \\ Universit\'a degli Studi 
di Milano--Bicocca\\ Via Cozzi, 53 20125 Milano \\ Italy}
\email{stefano.meda@unimib.it}

\author{Peter Sj\"ogren}
\address{ Department of Mathematical Sciences \\
G\"oteborg University and Chalmers \\ SE-412 96 G\"oteborg \\ Sweden}
\email{peters@math.chalmers.se}

\author{Maria Vallarino}
\address{Dipartimento di Matematica e Applicazioni\\Universit\'a degli Studi 
di Milano--Bicocca\\ Via Cozzi, 53 20125 Milano \\ Italy}
\email{maria.vallarino@unimib.it}

\thanks{Work partially supported by the
Progetto Cofinanziato ``Analisi Armonica''.}

\subjclass[2000]{42B30, 46A22}

\date{June 18, 2007 and, in revised form, October 05, 2007.}

\keywords{BMO, atomic Hardy space, extension of operators.}

\begin{abstract}
We prove that if $q$ is in $(1,\infty)$, $Y$ is a 
Banach space, and
$T$ is a linear operator defined on the space of finite linear combinations
of $(1,q)$-atoms in $\BR^n$ with the property that
$$
\sup\{ \norm{Ta}{Y}: \hbox{$a$ is a $(1,q)$-atom} \} < \infty,
$$
then $T$ admits a (unique) continuous extension to a bounded
linear operator from $\hu{\BR^n}$ to $Y$.  We show
that the same is true if we replace $(1,q)$-atoms by
\emph{continuous} $(1,\infty)$-atoms.  This is known
to be false for $(1,\infty)$-atoms.
\end{abstract}

\maketitle

\section{Introduction} \label{s: Introduction}

In a recent paper, M.~Bownik \cite{B2} showed that
there exists a linear functional $F$ defined on finite linear combinations
of $(1,\infty)$-atoms in $\BR^n$ with the property that
$$
\sup\{ \mod{F(a)}: \hbox{$a$ is a $(1,\infty)$-atom} \} < \infty,
$$
but which does not admit a continuous extension to $\hu{\BR^n}$.
If $v$ is a fixed function in $\lu{\BR^n}\setminus\{0\}$, then the operator
$B$, defined on finite linear combinations
of $(1,\infty)$-atoms by
$
Bf = F(f) \, v
$,
satisfies
$$
\sup\{ \norm{B a}{\lu{\BR^n}}: \hbox{$a$ is a $(1,\infty)$-atom} \} < \infty,
$$
but does not admit an extension to a bounded operator from $\hu{\BR^n}$
to $\lu{\BR^n}$.  This shows that the argument
``the operator $T$ maps $(1,\infty)$-atoms uniformly into $\lu{\BR^n}$,
and
hence it extends to a bounded operator from $\hu{\BR^n}$ to $\lu{\BR^n}$''
is fallacious.  

Fortunately, if $T$ is a Calder\'on--Zygmund operator, 
then the uniform boundedness of $T$ on 
$(1,\infty)$-atoms implies the boundedness from $\hu{\BR^n}$ to $\lu{\BR^n}$
(see, for instance, \cite[Ch.~7.3, Lemma~1]{MC},
\cite[Ch.~1.9]{B1}, \cite[Ch.~III.7]{GR}
and \cite[Thm 6.7.1]{G}).

The purpose of this paper is to show that the operator $B$ constructed
above is, to a certain extent, pathological.  Indeed,
we prove that if $q$ is in $(1,\infty)$, $Y$ is a Banach space, and
$T$ is a linear operator on finite linear combinations
of $(1,q)$-atoms in $\BR^n$ with the property that
\begin{equation}  \label{f: lim unif}
\sup\{ \norm{Ta}{Y}: \hbox{$a$ is a $(1,q)$-atom} \} < \infty,
\end{equation}
then $T$ admits a unique continuous extension to a bounded
linear operator from $\hu{\BR^n}$ to $Y$.
The same conclusion holds if we assume that
$T$ is a linear operator on finite linear combinations
of \emph{continuous} $(1,\infty)$-atoms in $\BR^n$ with the property that
\begin{equation}  \label{f: lim unif II}
\sup\{ \norm{Ta}{Y}: \hbox{$a$ is a continuous $(1,\infty)$-atom} \} 
< \infty.
\end{equation}
Note that this does not contradict Bownik's example.  Indeed, 
the restriction of the operator~$B$ to
continuous $(1,\infty)$-atoms extends
to a bounded operator $\wt B$ from $\hu{\BR^n}$ to~$\lu{\BR^n}$.
However, $B$ and $\wt B$ will agree on continuous $(1,\infty)$-atoms
but not on all $(1,\infty)$-atoms. 

To explain the idea of the proofs of these results, we need more notation.
Suppose that~$q$ is 
in $(1,\infty]$, and denote by $\huqfin{\BR^n}$
the vector space of all finite linear combinations of $(1,q)$-atoms.
Notice that $\huqfin{\BR^n}$ consists of all $\laq{\BR^n}$
functions with compact support and integral $0$.
Clearly, $\huqfin{\BR^n}$ is a dense subspace of $\hu{\BR^n}$.  
We may define a norm
on $\huqfin{\BR^n}$ as follows
$$
\norm{f}{\huqfin{\BR^n}}
= \inf \Bigl\{ \sum_{j=1}^N \mod{\la_j}: f = \sum_{j=1}^N \la_j \, a_j,
\, \, 
\hbox{$a_j$ is a $(1,q)$-atom, $N\in\BN$} \Bigr\}.
$$
Obviously $\norm{f}{\hu{\BR^n}} \leq  \norm{f}{\huqfin{\BR^n}}$
for every $f$ in $\huqfin{\BR^n}$.
An example due to Y.~Meyer, see \cite[p.~513]{MTW}, Bownik's paper \cite{B2} or
\cite[p.~370]{GR}, shows that
$\norm{\cdot}{\hu{\BR^n}}$ and $\norm{\cdot}{\huinftyfin{\BR^n}}$
are inequivalent norms on $\huinftyfin{\BR^n}$.  This is the 
starting-point of Bownik's construction. 

We prove that Meyer's example itself is somewhat exceptional.
Indeed, by using the maximal characterisation of $\hu{\BR^n}$,
we show that if $q<\infty$, then 
$\norm{\cdot}{\hu{\BR^n}}$ and $\norm{\cdot}{\huqfin{\BR^n}}$
are equivalent norms on $\huqfin{\BR^n}$ (see Section~\ref{s: Further}).  
Similarly, we prove that 
$\norm{\cdot}{\hu{\BR^n}}$ and $\norm{\cdot}{\huinftyfin{\BR^n}}$
are equivalent norms on $\huinftyfin{\BR^n}\cap C(\BR^n)$.  
This immediately
implies that operators defined on $\huqfin{\BR^n}$
which have either property (\ref{f: lim unif}) or
property (\ref{f: lim unif II}) automatically extend
to bounded operators from $\hu{\BR^n}$ to $\lu{\BR^n}$.

As discussed briefly in Section 3, this equivalence of norms remains true
for $H^p(\BR^n)$ with $0<p<1$ and $(p,q)$-atoms.

The extension property for operators was also proved, by different methods, for
$0 < p \leq 1$ and $(p,2)$-atoms and operators taking values in quasi-Banach
spaces, by D.~Yang and Y.~Zhou \cite{YZ}.


A theory of Hardy spaces has been developed in spaces of
homogeneous type; see R.R.~Coifman and G.~Weiss \cite{CW}.
 It is, however, not evident whether
our results extend to this case in general. Nevertheless,  let $M$
be such a space. By a simple functional analysis argument,
we show that if $q$ is in $(1,\infty)$
and $T$ is an operator defined on $\huqfin{M}$
satisfying the analogue of (\ref{f: lim unif}),
then $T$ automatically extends
to a bounded operator from $\hu{M}$ to $\lu{M}$ (see Section~\ref{s: Main
result}). 
It may be worth noticing that the proof of this
result applies also to certain metric measured
spaces $(M,\rho,\mu)$ where $\mu$ is only ``locally doubling''
\cite{MM, CMM, V}.

 For so-called RD-spaces, which are spaces of homogeneous type having  
``dimension $n$''
in a certain sense, our full results were recently extended
in the paper \cite{GLY} by L.~Grafakos, L.~Liu and Yang. These authors
consider $n/(n+1) < p \leq 1$ and quasi-Banach-valued operators.

The authors wish to thank N.Th.~Varopoulos for useful
conversations on the subject of this paper. 

\section{Notation and terminology} \label{s: Notation}

Suppose that $(M,\rho,\mu)$ is a space of homogeneous
type in the sense of Coifman and Weiss \cite{CW} and that
$\mu$ is a $\si$-finite measure.  
For the sake of simplicity, we shall assume that $\mu(M)$ is infinite.

Suppose that $q$ is in 
$(1,\infty]$.  For each closed ball $B$ in $M$,
we denote by $\laqO{B}$ the space of all functions in
$\laq{M}$ which are supported in $B$ and have integral~$0$.
Clearly $\laqO{B}$ is a closed subspace of $\laq{M}$.
The union of all spaces $\laqO{B}$ as $B$ varies
over all balls coincides with the space $\laqcO{M}$ of
all functions in $\laq{M}$ with compact support and integral~$0$.  
Fix a reference point $o$ in $M$ and for each positive integer $k$
denote by $B_k$ the ball centred at~$o$ with radius~$k$. 
A convenient way of topologising $\laqcO{M}$
is to interpret $\laqcO{M}$ as the strict inductive limit
of the spaces $\laqcO{B_k}$ (see \cite[II, p.~33]{Bo} for the definition of
the strict inductive limit topology).  
We denote by $X^q$ the space $\laqcO{M}$ with this topology,
and write~$X_k^q$ for $\laqcO{B_k}$.

We recall the basic definitions and results 
concerning the atomic Hardy space $\hu{M}$.
The reader is referred to \cite{CW} and the references therein
for this and more on Hardy spaces defined on spaces of homogeneous type.
Suppose that $q$ is in $(1,\infty]$.  A $(1,q)$-atom is a function $a$ 
in $\laq{M}$
supported in a ball $B$, with mean value~$0$ and such that 
$$
\Bigl(\frac{1}{\mu(B)} \, \int_B \mod{a}^q \wrt\mu\Bigr)^{1/q}
\leq \mu(B)^{-1}  
$$
if $q$ is finite, and   
$
\norm{a}{\infty}
\leq \mu(B)^{-1} 
$
if $q=\infty$.  We denote by $H^{1,q}(M)$ the space of all functions~$g$ 
in $\lu{M}$ which admit a decomposition of the form
$g=\sum_j \la_j\, a_j$, 
where the $a_j$ are $(1,q)$-atoms and the $\la_j$ are complex numbers 
such that $\sum_j \mod{\la_j} < \infty$. 
The norm $\norm{g}{H^{1,q}}$ of $g$ in $H^{1,q}(M)$ is 
the infimum of $\sum_j \mod{\la_j}$ over all such decompositions. 
It is well known that all the spaces $H^{1,q}(M)$
with $q\in (1,\infty)$ coincide with $H^{1,\infty}(M)$, 
and we denote them all by $\hu{M}$. 
Clearly,
the vector space $\huqfin{M}$ of all finite linear combinations
of $(1,q)$-atoms is dense in $\hu{M}$ with respect to 
the norm of $\hu{M}$, for $q$ in $(1,\infty]$. 
Observe also that $\huqfin{M}$ and $\laqcO{M}$ agree
as vector spaces, and so do the space of finite linear
combinations of continuous $(1,\infty)$-atoms and $\huinftyfin{M}\cap
C(\BR^n)$.

For each ball $B$ and each locally integrable function $f$,
we denote by $f_B$ the average of $f$ on~$B$. 
Recall that $BMO$ is the Banach space
of all locally integrable functions $f$, defined modulo
constants, such that 
$$
\norm{f}{BMO} 
= \sup_B \frac{1}{\mu(B)}\, \int_B \bigmod{f-f_B} \wrt \mu <\infty.
$$
The dual of $\hu{M}$ may be identified with $BMO$.

There are several characterisations of the space $\hu{\BR^n}$.
We shall make use of the so-called maximal characterisation,
which we briefly recall.  Suppose that $m$
is an integer with $m>n$, and denote by $\cA_m$ the
set of all functions $\vp$ in the Schwartz space $\cS(\BR^n)$
such that 
$$
\sup_{\mod{\be}\leq m}
\sup_{x\in \BR^n} (1+\mod{x})^m\, \bigmod{D^\be\vp(x)} \leq 1,
$$
where $\mod{\be}$ denotes the length of the multi-index $\be$.
For $\vp$ in $\cS(\BR^n)$ denote by $\vp_t$ the function
$t^{-n} \, \vp(\cdot/t)$.
Given $f$ in $\lu{\BR^n}$, define the ``grand maximal function''
$\cM_m f$ by
$$
\cM_mf
= \sup_{\vp \in \cA_m} \, \sup_{t>0} \mod{\vp_t\ast f}.
$$
The following result is classical \cite{FS, Sj, GR, St2}.

\begin{theorem}
Suppose that $f$ is in $\lu{\BR^n}$.  The following are equivalent:
\begin{enumerate}
\item[\itemno1]
$f$ is in $\hu{\BR^n}$;
\item[\itemno2]
the grand maximal function $\cM_mf$ is in $\lu{\BR^n}$.
\end{enumerate}
Furthermore, $f\mapsto \norm{\cM_m f}{\lu{\BR^n}}$ is an equivalent
norm on $\hu{\BR^n}$. 
\end{theorem}


The letter $C$ will denote a positive constant, 
which need not be the same at different occurrences.
Given two positive quantities $A$ and $B$, we shall mean by $A \sim B$ 
that there exists a constant $C$ such that $1/C\leq A/B\leq C$.


\section{The Euclidean case} \label{s: Further}

In this section we work in the classical setting of $\BR^n$.

\begin{theorem} \label{t: equivalence}
The following hold:
\begin{enumerate}
\item[\itemno1]
if $q<\infty$, then 
$\norm{\cdot}{\huqfin{\BR^n}}$ and $\norm{\cdot}{\hu{\BR^n}}$
are equivalent norms on $\huqfin{\BR^n}$; 
\item[\itemno2]
the two norms
$\norm{\cdot}{\huinftyfin{\BR^n}}$ and $\norm{\cdot}{\hu{\BR^n}}$
are equivalent on $\huinftyfin{\BR^n} \cap C(\BR^n)$. 
\end{enumerate}
\end{theorem}

\begin{proof}
Clearly, $\norm{f}{\hu{\BR^n}} \leq  \norm{f}{\huqfin{\BR^n}}$
for $f$ in $\huqfin{\BR^n}$ and for $q$ in $(1,\infty]$.
Thus, we have to show that for every $q$ in $(1,\infty)$
there exists a constant $C$ such that 
$$
\norm{f}{\huqfin{\BR^n}}
\leq  C \, \norm{f}{\hu{\BR^n}} 
\quant f \in \huqfin{\BR^n},
$$
and that a similar estimate holds for $q=\infty$
and all $f$ in $\huinftyfin{\BR^n} \cap C(\BR^n)$.

Suppose that $q$ is in $(1,\infty]$,
and that $f$ is in $\huqfin{\BR^n}$ with
$\norm{f}{\hu{\BR^n}} = 1$. By the translation
invariance of Lebesgue measure, we may assume that the support of
$f$ is contained in the closed ball $B=B(0,R)$ centred at $0$ with radius $R$.  
For each $k$ in $\BZ$, denote by $\Om_k$
the level set
$
\{x \in \BR^n: \cM_m f(x) > 2^k\}
$
of the grand maximal function $\cM_m f$ of $f$. 
We choose Whitney cubes $Q_i^k$, $i\in \BN$, with disjoint interiors 
satisfying $\Omega_k=\bigcup_iQ_i^k $ and
\begin{equation} \label{diamQ}
{\rm{diam}} (Q_i^k) \leq \eta \, {\rm{dist}}(Q_i^k, \Omega_k^c)\leq 4\, {\rm{diam}} (Q_i^k),
\end{equation}
where $\eta$ is a suitable constant in $(0,1)$.
Except for the factor $\eta$, this is Theorem VI.1 of \cite[p. 167]{St1}.
The only modification needed in the proof in \cite{St1} concerns the choice
of the constant denoted $c$.

By following closely  the proof  of \cite[Theorem~III.2~p.~107]{St2}
or \cite[Theorem~3.5~pp.~12-18]{Sj}, we produce an atomic decomposition of $f$ of the form
\begin{equation} \label{f: decomposition}
f = \sum_{i,k} \la_i^k \, a_i^k,
\end{equation}
such that the following hold:
\begin{enumerate}
\item[(a)]
$\mod{\la_i^k \, a_i^k} \leq C \, 2^k$ for every $k$ in $\BZ$;
\item[(b)]
for each $k$ in $\BZ$, the atoms $a_i^k$ are supported
in balls $B_i^k$ concentric with the $Q_i^k$ and contained in $\Om_k$.  
By choosing the constant $\eta$ in (\ref{diamQ}) small enough, depending
on the dimension, we can also ensure that
the family $\{B_i^k\}_{i}$ has the bounded overlap property,
uniformly with respect to $k$; 
\item[(c)]
there exists a constant $C$ independent of $f$ such that
$
\sum_{i,k} \bigmod{\la_i^k} \leq C \, \norm{f}{\hu{\BR^n}} = C.
$
\end{enumerate}
We write $2B$ for the closed ball concentric with $B$ whose radius is twice as
large. For $\vp$ in $\cA_m$
and $x$ in $\BR^n\setminus (2B)$ one then has
$$
\begin{aligned}
\mod{\vp_t\ast f(x)} 
& \leq t^{-n} \,\sup_{y\in B^c} \mod{\vp(y/t)} \, \norm{f}{\lu{\BR^n}} \\
& \leq t^{-n} \, (1+R/t)^{-m} \, \norm{f}{\lu{\BR^n}} \quant t \in \BR^+,
\end{aligned}
$$
so that
$$
\cM_mf(x)
= \sup_{\vp \in \cA_m} \, \sup_{t>R} \mod{\vp_t\ast f(x)} \\
\leq R^{-n},
$$
since $m>n$. Now, if $x$ is in $\Om_k \setminus(2B)$, the above inequality and the
definition of $\Om_k$ force
$2^k<R^{-n}$; denote by $k'$ the largest integer $k$
such that $2^k<R^{-n}$.  Then $\overline{\Om_k}$
is contained in $2B$ for $k>k'$.  

Next we define the functions $h$ and $\ell$ by
\begin{equation} \label{f: atomo A}
h = \sum_{k\leq k'} \sum_{i}  \la_i^k \, a_i^k
\qquad \hbox{and} \qquad
\ell =  \sum_{k> k'} \sum_{i} \la_i^k \, a_i^k. 
\end{equation}
Observe that both these series
converge in $\lu{\BR^n}$, simply because 
$\sum_{i,k} \bigmod{\la_i^k} < \infty$, so that~$h$
and $\ell$ have integral $0$.  
Clearly, $f = h + \ell$.
Furthermore, the support of
$\ell$ is contained in $2B$, because
it is contained in $\overline{\Om}_k$ by (b) above, and $\overline{\Om}_k$
is contained in $2B$ for all $k> k'$.
Therefore $h = f = 0$ in $(2B)^c$.

To estimate the size of $h$ in $2B$, we 
use (a) above and the bounded overlap property of~(b),
getting
$$
\mod{h}
 \leq C \sum_{k\leq k'} 2^k \; \leq \; C 2^{k'} \; \leq \; C  \mod{2B}^{-1}.
$$
This proves that $h/C$ is a $(1,\infty)$-atom,
where $C$ is independent of $f$.

Now we assume that $q<\infty$, and conclude the proof of \rmi.
Observe that $\ell$ is in $\laq{\BR^n}$, because $\ell = f - h$, and both
$f$ and $h$ are in $\laq{\BR^n}$.

We claim that the series
$\sum_{k>k'} \sum_i\la_i^k \, a_i^k$ converges to
$\ell$ in $\laq{\BR^n}$. 


Fixing $s$ in $\BZ$, we shall estimate 
$\sum_{k>k'} \sum_i \mod{\la_i^k \, a_i^k}$ 
in $\Om_{s}\setminus\Om_{s+1}$. First observe that all terms with $k>s$
vanish outside $\Om_{s+1}$. Then apply (a) and (b), 
to get the pointwise bound
$$
\sum_{k>k'} \sum_{i} \mod{\la_i^k \, a_i^k}
\; \leq \; C \sum_{k\leq s} 2^k 
 \;\leq  \; C 2^{s} 
 \;\leq \; C  \cM_mf.
$$
The constants $C$ above are independent
of $f$ and $s$, so that
$$
\sum_{k>k'} \sum_{i} \mod{\la_i^k \, a_i^k}
\leq C \, \cM_mf
$$
in all of $\BR^n$, with $C$ independent of $f$.  
Note that $\cM_mf$ is in $\laq{\BR^n}$, since $f$ is.  
This implies that the series defining $\ell$
converges almost everywhere, and the limit
must coincide with the $L^1$ limit~$\ell$.
The Lebesgue dominated convergence theorem now implies
that 
$\sum_{k>k'} \sum_{i} \la_i^k \, a_i^k$ 
converges to $\ell$ in $\laq{\BR^n}$, and
the claim is proved. 

Finally, for each positive integer $N$ we denote by $F_N$
the finite set of all pairs of integers $(i,k)$ such that $k>k'$
and $\mod{i}+\mod{k} \leq N$, and by
$\ell_N$ the function 
$
\sum_{(i,k) \in F_N} \la_i^k\, a_i^k.
$
The function $\ell_N$ is in $\huqfin{\BR^n}$,
and $f = h + \ell_N + (\ell-\ell_N)$.  Observe that $\ell-\ell_N$
will be a small multiple of a $(1,q)$-atom for large $N$.
Indeed, by taking $N$ large enough, we can make the corresponding 
coefficient less than any given $\vep$ in $\BR^+$.
Then
$$
\norm{f}{\huqfin{\BR^n}}
\; \leq \; C + \sum_{(i,k)\in F_N} \mod{\la_i^k} + \vep,
$$ 
so that 
$$
\norm{f}{\huqfin{\BR^n}}
\; \leq \; C + \sum_{(i,k)\in F_N} \mod{\la_i^k} 
\; \leq \; C, 
$$ 
by property (c) above, as required to conclude the proof of \rmi.

Now we finish the proof of \rmii.  
Assume that $f$ is a continuous function in $\huinftyfin{\BR^n}$.
A careful examination of the proof of \cite[Theorem III.2~pp.~107-8]{St2} 
or \cite[Theorem~3.5~pp.~12-18]{Sj}
shows that the atoms $a_i^k$ that appear in the decomposition 
(\ref{f: decomposition}) are then continuous. 
Furthermore,
we see that for each $k$ and $i$ the function $\la_i^k\, a_i^k$ 
depends only on the restriction of~$f$ to a ball $\tilde{B}_i^k$ 
which is a concentric enlargement of the ball $B_i^k$ from (b) above, 
by a fixed scaling factor. It is straightforward to check that
if $f$ is constant in $\tilde{B}_i^k$, then~$\la_i^k\, a_i^k=0$, and
that there exists an absolute constant~$C$ such that 
if $\mod{f} < \vep$ in $\tilde{B}_i^k$, then 
$\bigmod{\la_i^k\, a_i^k} < C \, \vep$.

Since trivially $\cM_mf \leq C_n\norm{f}{\infty}$ where the constant
$C_n$ depends only on $n$,
the level set $\Om_k$ is empty for all $k$ such that 
$2^k\geq C_n \norm{f}{\infty}$.  We denote by $k''$ the largest integer
for which the last inequality does not hold. Then the index $k$ in 
the sum defining $\ell$ in (\ref{f: atomo A}) will run only over 
$k' < k \leq k''$. 

Let $\vep$ be positive. Since $f$ is uniformly continuous, there
exists a positive $\delta$ such that $|x-y| < \delta$ implies
$$
\bigmod{f(x) - f(y)} < \vep.
$$
Write $\ell = \ell_1^\vep+\ell_2^\vep$ with
$$
\ell_1^\vep = \sum_{(i,k) \in F_1} \la_i^k \, a_i^k
\qquad\hbox{and} \qquad
\ell_2^\vep = \sum_{(i,k) \in F_2} \la_i^k \, a_i^k,
$$
where 
$
F_1 = \{(i,k): \diam(\tilde{B}_i^k) \geq \de,\; k' < k \leq k''\}
$
and
$
F_2 = \{(i,k): \diam(\tilde{B}_i^k) < \de,\; k' < k \leq k''\}.
$
Since $F_1$ is a finite set,  $\ell_1^\vep$
is  continuous.

To estimate $\ell_2^\vep$, we  denote by $x_i^k$ the centre of the ball $B_i^k$ 
and write for $(i,k)$ in $F_2$
$$
f(x) 
= f(x_i^k) + f(x) - f(x_i^k). 
$$
Then $\mod{\la_i^k\, a_i^k} < C \, \vep$,
because $\bigmod{f(x) - f(x_i^k)} < \vep$ for $x$ in $\tilde{B}_i^k$. 
For fixed $k$ the balls $\{B_i^k\}_{i}$
have  uniformly bounded overlap, so there exists an absolute constant 
$C$ such that
$$
\mod{\ell_2^\vep}
\leq C \sum_{k' < k \leq k''} \vep
\leq C \, (k''-k') \, \vep.
$$
Since $\vep$ is arbitrary, we can thus split $\ell$ into a continuous part
and a part that is uniformly arbitarily small. It follows that  $\ell$
is continuous. But then $h=f-\ell$ is also continuous, so that $h$
 is  a continuous $(1,\infty)$-atom, multiplied by a factor $C$.

To find a finite atomic decomposion of $\ell$, we use again the splitting
$\ell = \ell_1^\vep+\ell_2^\vep$. Clearly $\ell_1^\vep$ is for each $\vep$ 
a finite linear combination of continuous $(1,\infty)$-atoms, and the
$\ell^1$ norm of the coefficients is controlled by $\|f\|_{H^1}$, in 
view of (c). Observe that $\ell_2^\vep = \ell - \ell_1^\vep$ is continuous.
Further, $\ell_2^\vep$ is supported in $2B$, has integral 0 and satisfies
$|\ell_2^\vep| \leq C(k''-k')\vep$. Choosing $\vep$, we can thus make
 $\ell_2^\vep$ into an arbitrarily small multiple of a continuous 
$(1,\infty)$-atom.

To sum up, $f = h + \ell_1^\vep + \ell_2^\vep$ gives the desired finite
atomic decomposition of $f$, with coefficients controlled by 
$\|f\|_{H^1}$. 

We have completed the proof of (ii) and that of the theorem. 
\end{proof}

\begin{remark}
Theorem \ref{t: equivalence} (ii) implies that any function $f$ 
in $\huinftyfin{\BR^n} \cap C(\BR^n)$ admits a finite 
decomposition in $(1,\infty)$-atoms such that the sum of 
the corresponding coefficients is 
$\leq C\,\norm{f}{\hu{\BR^n}}$. Actually, the proof of 
Theorem \ref{t: equivalence} (ii) shows that 
we can construct this finite decomposition in such a way that it 
involves only continuous $(1,\infty)$-atoms.
\end{remark}

\begin{remark}
Theorem  \ref{t: equivalence} extends to $H^p(\BR^n)$ with $0<p<1$ and $(p,q)$-atoms, 
where one
can now have $1\leq q \leq \infty$. The proof is rather similar  to the one given
above, so we  only describe briefly the  modifications needed for
part (i). Let thus $1 \leq q < \infty$.
Given $f \in H^{p,q}_{\mathrm{fin}}(\BR^n)$ supported in a ball $B_R$, the first step
is the inequality $\mathcal{M}_mf \leq CR^{-n/p} \|f\|_{H^p(\BR^n)}$, valid outside a 
larger ball $B_{CR}$. One proves this by comparing the values of 
$\mathcal{M}_mf$
at different points and using the fact that 
$\|\mathcal{M}_mf\|_{L^p(\BR^n)}\sim \|f\|_{H^p(\BR^n)}$. 
Then the $\Omega_k$ and the decompositions $f = \sum \lambda_i^k a_i^k = h+\ell$ are 
introduced as above. The sum $\ell$ now converges in $\mathcal{S}'$ and 
is dominated by $\mathcal{M}_mf$. If $q>1$, we have 
$\mathcal{M}_mf \in L^q(\BR^n)$ and conclude
as before that  $\ell$  converges in $L^q(\BR^n)$. For $q=1$, the tail sum
$S_\kappa = \sum_{k\geq \kappa} \sum_i \lambda_i^k a_i^k$ tends to 0 in $L^1(\BR^n)$
as $\kappa \to +\infty$,
because $S_\kappa$ is nonzero only in $\Omega_\kappa$ and not larger than
$|f| + C2^\kappa$ there, and $|\Omega_\kappa| = o(2^{-\kappa})$
as $\kappa \to +\infty$. The rest of the proof goes as before.
See also \cite[Theorem 5.6]{GLY}
\end{remark}

\begin{corollary} \label{c: operators}
Suppose  that $Y$ is a Banach space and that one of the following holds:
\begin{enumerate}
\item[\itemno1]
$q$ is in $(1, \infty)$ and $T: \huqfin{\BR^n} \to Y$
is a linear operator such that
$$
A : = \sup\{ \norm{T a}{Y}: \hbox{$a$ is a $(1,q)$-atom} \} < \infty;
$$
\item[\itemno2]
$T$ is a $Y$-valued linear operator defined on continuous $(1,\infty)$-atoms 
such that
$$
A : = \sup\{ \norm{T a}{Y}: \hbox{$a$ is a continuous
$(1,\infty)$-atom} \} < \infty.
$$
\end{enumerate}
Then there exists a unique bounded linear operator
$\wt T$ from $\hu{\BR^n}$ to $Y$ which extends~$T$.
\end{corollary}

\begin{proof}
We consider the case \rmi.
Suppose that $f$ is in $\huqfin{\BR^n}$, $f = \sum_{j=1}^N 
\la_j \, a_j$ say, where $a_j$ are $(1,q)$-atoms.  Then
the assumption and the triangle inequality give 
$$
\norm{Tf}{Y}
\leq A \, \sum_{j=1}^N \mod{\la_j}.
$$
By taking the infimum of the right-hand side with respect
to all decomposition of $f$ as a finite sum of $(1,q)$-atoms,
we obtain
$$
\norm{Tf}{Y}
\leq A \, \norm{f}{\huqfin{\BR^n}}.
$$
Now, Theorem~\ref{t: equivalence}~\rmi\ implies that the right-hand
side is dominated by $CA\, \norm{f}{\hu{\BR^n}}$,
where~$C$ does not depend on $f$, 
and a density argument completes the proof of the corollary.

The case \rmii\ is similar.
\end{proof}

\begin{remark}
The statement of Corollary~\ref{c: operators}~\rmi\
becomes false if we replace $q$ by $\infty$.  A counterexample
is given by the operator $B$ defined in the Introduction.
Note also that Corollary~\ref{c: operators} applies
to linear functionals.
\end{remark}

\section{Results on spaces of homogeneous type} \label{s: Main result}

In this section, we work in a space of homogeneous type
$(M,\rho,\mu)$.  Recall that we assume that $\mu$
is $\si$-finite, and that $\mu(M)$ is infinite. 

\begin{theorem} \label{t: basic prop}
Suppose that $q$ is in $(1,\infty)$ and that
$T$ is a linear operator defined on $\huqfin{M}$ with the
property that
$$
A : = \sup\{ \norm{T a}{\lu{M}}: \hbox{$a$ is a $(1,q)$-atom} \} < \infty.
$$
Then there exists a unique bounded linear operator
$\wt T$ from $\hu{M}$ to $\lu{M}$ which extends~$T$.
\end{theorem}

\begin{proof}
We prove the result in the case where $q=2$.  The proof 
in the other cases is similar.

Suppose that $B$ is a ball.  For each
$f$ in $\ldO{B}$ such that $\norm{f}{\ld{M}}=1$,
the function
$
\mu (B)^{-1/2} \, f 
$
is a $(1,2)$-atom, so that 
$$
\norm{Tf}{\lu{M}}
\leq A\, \mu (B)^{1/2} \quant f \in \ldO{B}
$$
by the assumption.
In particular, the restriction of $T$ to $X_k^2$ is bounded from
$X_k^2$ to $\lu{M}$ for each $k$.  Thus, $T$ is bounded from
$X^2$ to $\lu{M}$. It follows that $T^*$ is
bounded from $\ly{M}$ to the dual of $X^2$.
But the dual of $X^2$ is the quotient space $\ldloc{M}/\BC$, since 
that of  $L^2_{c,0}(B_k)$  is  $L^2(B_k)/\BC$. 
Now, for every $f$ in $\ly{M}$ and for every $(1,2)$-atom $a$
$$
\prodo{T a}{f}
= \prodo{a}{T^* f} \\
= \int_M {a}\, {T^* f} \wrt \mu,
$$
so that
$$
\Bigmod{\int_M {a}\, {T^* f} \wrt \mu}
 = \bigmod{\prodo{T a}{f}} 
 \leq A \, \norm{f}{\infty}.
$$
A standard argument then shows that $T^* f$ belongs to $BMO(M)$ and that
\begin{equation}  \label{f: BMO I}
\norm{T^* f}{BMO(M)}
\leq 2A \, \norm{f}{\infty}
\quant f \in \ly{M}.
\end{equation}
We give the details for the reader's convenience.
Suppose that $B$ is a ball and observe that
$$
\Bigl[\int_B \mod{T^*f-(T^*f)_B}^2 \wrt \mu\Bigr]^{1/2}
= \sup_{\norm{\vp}{\ld{B}} = 1}
      \Bigmod{\int_B \vp \, \bigl(T^*f-(T^*f)_B\bigr) \wrt \mu}.
$$
But
$$
\begin{aligned}
\int_B \vp \, \bigl(T^*f-(T^*f)_B\bigr) \wrt \mu
& = \int_B \bigl(\vp-\vp_B\bigr) \, \bigl(T^*f-(T^*f)_B\bigr) \wrt \mu \\
& = \int_B \bigl(\vp-\vp_B\bigr) \, T^*f \wrt \mu,
\end{aligned}
$$
and since $\norm{\vp}{\ld{B}} = 1$
$$
\bigmod{\vp_B}
 \leq \Bigl[\frac{1}{\mu(B)} \, \int_B \mod{\vp}^2 \wrt \mu\Bigr]^{1/2} \;
 \leq \;\mu(B)^{-1/2}.
$$
Write $\psi$ instead of $\vp-\vp_B$.  Then 
$$
\norm{\psi}{\ld{B}}
\leq \norm{\vp}{\ld{B}} + \mod{\vp_B} \, \mu(B)^{1/2} \\
\leq 2,
$$
so that $\psi/(2\,\mu(B)^{1/2})$ is a $(1,2)$-atom. 
Therefore
$$
\Bigmod{\int_B \psi \,\,  T^*f \wrt \mu}
\leq 2A\, \mu(B)^{1/2} \, \norm{f}{\infty}
$$
Combining the above, we conclude that
for every ball $B$
$$
\Bigl[\frac{1}{\mu(B)} \, 
\int_B \mod{T^*f-(T^*f)_B}^2 \wrt \mu\Bigr]^{1/2}
\leq 2A\, \norm{f}{\infty},
$$
and (\ref{f: BMO I}) follows.

Now we show that $T$ extends
to a bounded operator from $\hu{M}$ to $\lu{M}$
with norm at most $2A$.
Observe that $X^2$ and $\hudfin{M}$ coincide as vector spaces.
For every $g$ in $\hudfin{M}$ and for every $f$ in $\ly{M}$
$$
\begin{aligned}
\bigmod{\prodo{Tg}{f}}
& =    \bigmod{\prodo{g}{T^*f}} \\
& \leq \norm{g}{\hu{M}} \, \norm{T^*f}{BMO(M)} \\
& \leq 2A\, \norm{g}{\hu{M}} \, \norm{f}{\ly{M}}. 
\end{aligned}
$$
By taking the supremum of both sides
over all functions $f$ in $\ly{M}$ with $\norm{f}{\ly{M}}=1$, we obtain that
$$
\norm{Tg}{\lu{M}}
\leq 2A\, \norm{g}{\hu{M}} 
\quant g \in \hudfin{M}.
$$
Finally we observe that $\hudfin{M}$ is dense in $\hu{M}$
(with respect to the norm of $\hu{M}$), and the required conclusion
follows by a density argument.
\end{proof}

Quite often one encounters the following situation.
Suppose that $T$ is a bounded linear operator on $\ld{M}$. 
Then $T$ is automatically
defined on 
$\hudfin{M}$. Assume that
$$
A : = \sup\{ \norm{T a}{\lu{M}}: \hbox{$a$ is a $(1,2)$-atom} \} < \infty.
$$
By the previous result, the restriction of $T$ to $\hudfin{M}$
has a unique extension to a bounded linear operator $\wt T$
from $\hu{M}$ to $\lu{M}$. The question is whether the operators
$T$ and $\wt T$ are consistent, i.e., whether they coincide on the 
intersection 
$\hu{M} \cap \ld{M}$ of their domains.
The answer to this question is in the affirmative as the following
proposition shows.

\begin{proposition}
Suppose that $T$ is bounded on $\ld{M}$ and that
$$
A : = \sup\{ \norm{T a}{\lu{M}}: \hbox{$a$ is a $(1,2)$-atom} \} < \infty.
$$
Denote by $\wt T$ the unique continuous linear extension
of the restriction of \ $T$ to $\hudfin{M}$
to an operator from $\hu{M}$ to $\lu{M}$.
Then the operators $T$ and $\wt T$ agree on $\hu{M} \cap \ld{M}$.
\end{proposition}

\begin{proof}
Suppose that $f$ is in $\ld{M} \cap \ly{M}$ and that
$g$ is in $\ldcO{M}$. Denote by $T^*$ the transpose operator
of $T$ (as an operator on $\ld{M}$). Then
\begin{equation} \label{f: prima}
\int_M {g}\,\, {T^*f} \wrt \mu
= \int_M {Tg}\,\, {f} \wrt\mu.
\end{equation}
Since $g$ is in $\hudfin{M}$ and the operators $T$ and $\wt T$
agree on $\hudfin{M}$, we see that
\begin{equation} \label{f: seconda}
\begin{aligned}
\int_M {Tg}\,\, {f} \wrt\mu
& = \int_M {\wt Tg}\,\, {f} \wrt\mu \\
& = \prodo{g}{(\wt T)^*f},
\end{aligned}
\end{equation}
where $(\wt T)^*$ denotes the transpose of the operator $\wt T$
from $\hu{M}$ to $\lu{M}$.
Note that $(\wt T)^*f$ is in $BMO(M)$ and $g$ is
a multiple of an atom.
Thus the above scalar product  $\prodo{g}{(\wt T)^*f}$ (with respect to the 
duality between $\hu{M}$
and $BMO(M)$) may be written as $\int_M {g}\, \,{(\wt T)^*f} \wrt\mu$.
Therefore, (\ref{f: prima}) and (\ref{f: seconda}) imply that
$$
\int_M {g}\, \,\bigl[T^*f - {(\wt T)^*f}\bigr] \wrt\mu
= 0
\quant g \in \ldcO{M},
$$
i.e., for all $g$ in $X^2$. Therefore
$T^*f - {(\wt T)^*f} = 0$ in the dual space of
$X^2$, i.e., in $\ldloc{M}/\BC$.
This implies that $T^*f - {(\wt T)^*f}$ is constant.

Now, suppose that $g$ is in $\hu{M}\cap \ld{M}$ and that
$f$ is in $\ld{M}\cap\ly{M}$.
\begin{equation} \label{f: terza}
\begin{aligned}
\int_M {Tg}\,\, {f} \wrt\mu
& = \int_M {g}\, \,{T^*f} \wrt\mu \\
& = \int_M {g}\, \,{(\wt T)^*f} \wrt\mu \\
& = \int_M {\wt Tg}\, \,{f} \wrt\mu.
\end{aligned}
\end{equation}
Since $f$ is an arbitrary function in $\ld{M}\cap\ly{M}$,
$Tg-\wt Tg = 0$ almost everywhere, as required.
\end{proof}

\bibliographystyle{amsplain}

\end{document}